\documentclass[11pt]{article}
\usepackage{amsmath, amssymb, cite}
\usepackage{graphicx, ifpdf, colordvi}

\newtheorem{thm}{Theorem}

\makeatletter
\def\ExtendSymbol#1#2#3#4#5{\ext@arrow 0099{\arrowfill@#1#2#3}{#4}{#5}}
\makeatother

\title{\bf\LARGE An Iterated Map for the Lebesgue Identity}
\author{William Y.C. Chen\raisebox{5pt}{\scriptsize 1},
Qing-Hu Hou\raisebox{5pt}{\scriptsize 2}, and Lisa H. Sun\raisebox{5pt}{\scriptsize 3}}
\date{Center for Combinatorics, LPMC-TJKLC\\
 Nankai University\\
Tianjin 300071, P.R. China \\
\vspace{15pt}
\raisebox{5pt}{\scriptsize 1\,}chen@nankai.edu.cn,
\raisebox{5pt}{\scriptsize 2\,}hou@nankai.edu.cn, \raisebox{5pt}{\scriptsize 3\,}sunhui@nankai.edu.cn}

\begin{document}
\maketitle

\noindent {\bf Abstract.} We present a simple iteration for the Lebesgue identity
                              on partitions, which leads to a refinement involving the alternating sums of partitions.



\vskip 6mm

We find a simple iterated map for the classical Lebesgue identity on partitions.
As an application of this iterated map, we give a refinement of
the partition interpretation of this identity involving alternating
sums of partitions. Recall that the $q$-shifted factorials
are defined by
\begin{equation*}
(a;q)_\infty=\prod_{k=0}^\infty (1-aq^k)\quad \mbox{and} \quad
(a;q)_{n}=\frac{(a;q)_\infty}{(aq^n;q)_\infty},\quad n\in
\mathbb{Z},
\end{equation*}
where $|q|<1$.
The Lebesgue identity reads
\begin{equation}\label{Lebes}
\sum_{k=0}^\infty \frac{(-aq;q)_k}{(q;q)_k} q^{k+1 \choose 2} = (-aq^2;q^2)_\infty (-q;q)_\infty,
\end{equation}
see, for example, Andrews \cite{Andrews98}.
There are several  combinatorial proofs of the Lebesgue identity.   Ramamani and Venkatachaliengar \cite{RamVen72} found a bijection for the following generalization
of (\ref{Lebes}),
\[
\sum_{m=0}^\infty q^{m(m+1)/2} \frac{(z;q)_m}{(q;q)_m} \alpha^m=(z;q)_\infty (-\alpha q;q)_\infty \sum_{n=0}^\infty \frac{z^n}{(q;q)_n(-\alpha q;q)_n}.
\]
Bessenrodt \cite{Bessenrodt94} gave a combinatorial interpretation
 in terms of $2$-modular diagrams.   Alladi and Gordon \cite{AllGor93}
 provided another bijection which implies the Lebesgue identity.
  Pak modified the construction of Alladi and Gordon to give
   a direct correspondence by using standard MacMahon diagrams \cite{Pak06}.
 Fu \cite{Fu08} discovered a   bijective proof of the following
 extension of (\ref{Lebes}) by applying the insertion algorithm of Zeilberger:
 \begin{equation*}
 \sum_{n=0}^\infty \frac{(-aq;q)_n}{(q;q)_n} b^n q^{n+1\choose 2}=(-bq;q)_\infty \sum_{k=0}^\infty \frac{(ab)^kq^{k(k+1)}}{(q;q)_k (-bq;q)_k}.
 \end{equation*}
Rowell \cite{Rowell09} presented a combinatorial proof which leads to the
 following finite form of (\ref{Lebes}):
 \begin{equation*}
 \sum_{n=0}^L {L\brack n}_q (-aq;q)_n q^{n(n+1)/2}=\sum_{k=0}^L {L \brack k}_{q^2} (-q;q)_{L-k} a^k q^{k(k+1)}.
 \end{equation*}  Recently, Little and Sellers \cite{LitSel09} have established
 the relation
\eqref{Lebes} by using weighted Pell tilings.

To describe our bijection, we follow the terminology  in \cite{Andrews98}.
A partition is meant to be a non-increasing finite sequence of positive integers $\lambda=(\lambda_1,\ldots,\lambda_\ell)$. The entries $\lambda_i$ are called the  parts of $\lambda$. The number of parts of $\lambda$ is denoted by $\ell(\lambda)$, and sum of parts is denoted by $|\lambda|=\lambda_1+\cdots + \lambda_\ell$.  The conjugate partition of $\lambda$ is denoted by $\lambda'$. The  partition with no parts is denoted by $\varnothing$.

\begin{figure}[ht]
\centering
  \setlength{\unitlength}{0.05 mm}%
  \begin{picture}(980.0, 583.1)(0,0)
  \put(0,0){\includegraphics{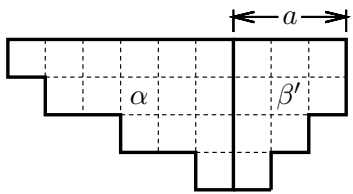}}
  \put(365.00,246.40){\fontsize{11.38}{13.66}\selectfont \makebox(40.0, 80.0)[l]{$\alpha$\strut}}
  \put(754.00,246.40){\fontsize{11.38}{13.66}\selectfont \makebox(40.0, 80.0)[l]{$\beta'$\strut}}
  \put(770.28,460.71){\fontsize{11.38}{13.66}\selectfont \makebox(120.0, 80.0)[l]{$a$\strut}}
  \end{picture}%
\caption{\label{lebesgue}A pair $(\alpha, \beta) \in P$}
\end{figure}

Denote the left hand side of the Lebesgue identity \eqref{Lebes} by $f(a,q)$. It is easily seen that
\[
f(a,q) = \sum_{(\alpha, \beta) \in P} a^{\ell(\beta)} q^{|\alpha|+|\beta|},
\]
where
$P$ denotes the set of pairs $(\alpha, \beta)$ of partitions with distinct parts
 such that  $\ell(\alpha)$ is not less than the largest part of $\beta$. The corresponding diagram is illustrated by Figure~\ref{lebesgue}.

Clearly, the right hand side of \eqref{Lebes} has the following combinatorial interpretation
\[
\sum_{(\mu,\nu) \in Q} a^{\ell(\nu)} q^{|\mu|+|\nu|},
\]
where $Q$ is the set of pairs $(\mu, \nu)$ of partitions with distinct parts
such that  $\nu$  has only even parts.

\begin{figure}[ht]
\centering
  \setlength{\unitlength}{0.05 mm}%
  \begin{picture}(2436.5, 1070.5)(0,0)
  \put(0,0){\includegraphics{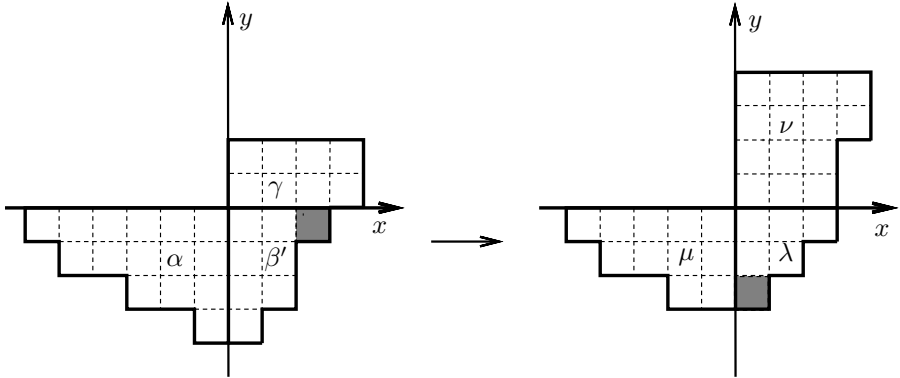}}
  \put(470.00,315.76){\fontsize{10.24}{12.29}\selectfont \makebox(36.0, 72.0)[l]{$\alpha$\strut}}
  \put(730.00,315.76){\fontsize{10.24}{12.29}\selectfont \makebox(36.0, 72.0)[l]{$\beta'$\strut}}
  \put(740.16,510.70){\fontsize{10.24}{12.29}\selectfont \makebox(36.0, 72.0)[l]{$\gamma$\strut}}
  \put(2099.40,670.49){\fontsize{10.24}{12.29}\selectfont \makebox(36.0, 72.0)[l]{$\nu$\strut}}
  \put(1830.70,330.37){\fontsize{10.24}{12.29}\selectfont \makebox(36.0, 72.0)[l]{$\mu$\strut}}
  \put(2095.46,330.06){\fontsize{10.24}{12.29}\selectfont \makebox(36.0, 72.0)[l]{$\lambda$\strut}}
  \put(1016.47,402.66){\fontsize{10.24}{12.29}\selectfont \makebox(36.0, 72.0)[l]{$x$\strut}}
  \put(663.30,956.34){\fontsize{10.24}{12.29}\selectfont \makebox(36.0, 72.0)[l]{$y$\strut}}
  \put(2352.46,396.42){\fontsize{10.24}{12.29}\selectfont \makebox(36.0, 72.0)[l]{$x$\strut}}
  \put(2016.77,952.43){\fontsize{10.24}{12.29}\selectfont \makebox(36.0, 72.0)[l]{$y$\strut}}
  \end{picture}%
\caption{\label{lebcase1} An example}
\end{figure}

For a triple of partitions $(\alpha, \beta, \gamma)$ where $(\alpha, \beta) \in P$ and $\gamma$ is a partition with even parts such that $\ell(\gamma)\geq \ell(\beta)$ or $\gamma=\phi$, we define a map $\phi \colon (\alpha, \beta, \gamma)\rightarrow (\mu, \lambda, \nu)$ as follows:

Case 1: The smallest part of $\beta$ equals $1$.
       Decrease each part of $\alpha$ by $1$ to
            form a partition $\mu$. Change  the $1$-part of $\beta$ to an $(\ell(\alpha)+1)$-part and  decrease each part of the resulting
             partition by $2$ to generate a partition $\lambda$. Then add two $\ell(\beta)$-parts to the conjugate partition $\gamma'$ to produce  a conjugate partition $\nu'$. This
              operation can be visualized as moving the rightmost square of $\beta'$ to the bottom of $\alpha$,
              then shifting the diagram below the $x$-axis to the right by one column, and finally moving up the diagram on the right side of the $y$-axis by two rows. See Figure \ref{lebcase1} for an illustration, where
             $\alpha=(6,5,3,1)$, $\beta=(4,3,1)$, $\gamma=(2,2,2,2)$, $\mu=(5,4,2)$, $\lambda=(3,2,1)$, and $\nu=(4,4,4,2)$.

Case 2: The smallest part of $\beta$ is larger than $1$.
 Set $\mu=\alpha$ and move up the diagram of $\beta'$ by two rows to generate the two conjugate partitions $\lambda'$ and $\nu'$.

To recover $(\alpha,\beta,\gamma)$ from $(\mu,\lambda,\nu)$, we first move down the diagram on the right side of the $y$-axis by two rows
to obtain a triple $(\bar{\alpha}, \bar{\beta}, \bar{\gamma})$.
If $\bar{\beta}_1 \leq \ell(\bar{\alpha})$, we then have $(\alpha,\beta,\gamma)=(\bar{\alpha}, \bar{\beta}, \bar{\gamma})$.
Otherwise, we further shift the diagram below the $x$-axis to the left by one column and move the bottom square to the right of $\bar{\beta}'$.
Thus, $\phi$ is invertible.

Starting from $(\alpha, \beta, \varnothing)$, we can iterate the above map until $\lambda$ becomes empty.  This gives a pair $(\mu, \nu)$ of  partitions
that  belongs to $Q$.
This completes the combinatorial proof of the Lebesgue identity.

The above map leads to a refinement of the Lebesgue identity \eqref{Lebes}.
Define the alternating sum of a partition $\lambda$ by
\[
|\lambda|_a=\lambda_1-\lambda_2+\lambda_3-\lambda_4+\cdots.
\]
This statistic has occurred in the study of  refinements of Euler's
partition theorem, see \cite{BousErik97}. Notice that when the parts of $\alpha$ are distinct, the alternating sum of $\alpha$ equals to the number of odd parts of its conjugate partition. Denote by $n_o(\lambda)$ the number of odd parts of a partition $\lambda$. It is straightforward to check that the map $\phi$ preserves the difference $n_o(\alpha')-n_o(\beta)$.
Therefore, our bijection leads to the following refinement of the combinatorial interpretation of the Lebesgue identity.

\begin{thm}\label{prop1}
Let $P$ denote the set of pairs $(\alpha, \beta)$ of partitions with distinct parts
such that  $\ell(\alpha)$ is not less than the largest part of $\beta$, and let $Q$ denote the set of pairs $(\mu, \nu)$ of partitions with distinct parts such that  $\nu$  has only even parts. Then for each nonnegative integer $k$, the number of pairs $(\alpha, \beta)\in P$ with $|\alpha|_a-n_o(\beta)=k$ is equal to the number of pairs $(\mu,\nu)\in Q$ with $|\mu|_a=k$.
\end{thm}

We notice that Bessenrodt's bijection \cite{Bessenrodt94} also keeps the difference $|\alpha|_a - n_o(\beta)$. We also note that our map can be viewed
as a direct correspondence in the sense that it does not
require Sylvester's bijection for Euler's identity, see the remark in \cite{Pak06}.

\vskip 15pt \noindent {\small {\bf Acknowledgments.}  We thank Byungchan Kim for helpful comments. This
work was supported by the National Science Foundation, the PCSIRT project
of the Ministry of Education, and the Ministry of Science and
Technology of China.}

\end{document}